\documentclass[10pt]{article}

\usepackage[superscript,sort]{cite}

\usepackage{ifthen,ifpdf}
\ifpdf
	\usepackage[pdftex]{hyperref}	 \hypersetup{colorlinks=true,linkcolor=black,citecolor=black,urlcolor=blue,backref=page,bookmarks=true,breaklinks=true,plainpages=false }
	\usepackage[pdftex]{graphicx}
	\usepackage[usenames,pdftex]{color}
\else
	 \usepackage[colorlinks=true,breaklinks=true,plainpages=false]{hyperref}
	\usepackage{graphicx}
	\usepackage[usenames]{color}
\fi

\usepackage{amsthm,amscd,amsxtra,amsfonts,amsmath,amssymb,multirow}
\usepackage{wrapfig}
\usepackage[footnotesize]{caption}
\usepackage[tiny,compact]{titlesec}
\usepackage[textwidth=0.8in,textsize=footnotesize]{todonotes}

\usepackage{helvet}

\usepackage[letterpaper,top=0.5in,bottom=0.5in,left=0.5in,right=0.5in]{geometry}

\titlespacing*{\section}{0pt}{*0}{*0}
\titlespacing*{\subsection}{0pt}{*0}{*0}
\titlespacing*{\subsubsection}{0pt}{*0}{*0}
\titlespacing{\paragraph}{0pt}{*0}{*1}
\definecolor{MyPurple}{rgb}{1,0,1}

\newtheorem{lemma}{Lemma}[section]

\begin{document}
\pagenumbering{roman}




\title{Accurate,  robust and  reliable  calculations of Poisson-Boltzmann solvation   energies}

\author{Bao Wang$^{1}$, and Guo-Wei Wei$^{1,2,3}$
\footnote{
Address correspondences  to Guo-Wei Wei. E-mail: wei@math.msu.edu}  \\
$^1$Department of Mathematics \\
Michigan State University, MI 48824, USA\\
$^2$Department of Electrical and Computer Engineering \\
Michigan State University, MI 48824, USA \\
$^3$Department of Biochemistry and Molecular Biology\\
Michigan State University, MI 48824, USA
}

\date{\today}

\maketitle

\begin{abstract}
Developing accurate solvers for the Poisson Boltzmann (PB) model is the first step to make the PB model suitable for implicit solvent simulation. Reducing the grid size influence on the performance of the solver benefits to increasing the speed of solver and providing accurate electrostatics analysis for solvated molecules. In this work, we explore the accurate coarse grid PB solver based on the Green's function treatment of the singular charges, matched interface and boundary (MIB) method for treating the geometric singularities, and posterior electrostatic potential field extension for calculating the reaction field energy. We made our previous PB software, MIBPB, robust and provides almost grid size independent reaction field energy calculation. Large amount of the numerical tests verify the grid size independence merit of the MIBPB software. The advantage of MIBPB software directly make the acceleration of the PB solver from the numerical algorithm instead of utilization of advanced computer architectures. Furthermore, the presented MIBPB software is provided as a free online sever.
\end{abstract}

\newpage

\label{introduction}
\section{Introduction}
Accurate and efficient modeling and computation of the electrostatics interaction of the biomolecule in solvent environment is important. First, it is crucial for studying many cellular process, such as signal transduction, gene expression and protein synthesis \cite{Harris:2013}. Second, it plays a key role in determining the structure and activity of the biomolecules \cite{Warwicker:1982, Warshel:1991, Warshel:1998}. Third, it has many pharmaceutical applications \cite{Harris:2013}, in the computer aided drug design, binding energy calculation relies on the accurate calculation of the electrostatics interaction of the biomolecule in the solvent environment. Fourth, it enables more realistic aqueous environment molecule dynamics simulation \cite{Davis:1990a}. Many methodologies have been developed in the past decades for modeling the electrostatics interaction in the solvent environment, these models are typically classified into two categories, the explicit solvent model and implicit solvent model. Both solute and solvent are modeled with atomic details in the explicit solvent model. The explicit solvent model is usually accurate, but it is extremely time consuming\cite{LuoRayPB-VS-Explicitsolvent:2006, LuoRayPB-VS-Explicitsolvent:2010}.
The implicit solvent model which is generally regarded as the best compromise between the accuracy and efficiency, which models the solute molecule at the atomic level, while solvent is modeled as a dielectric continuum. The implicit solvent model provides an effective approach for modeling the solvation phenomena, especially powerful in modeling the electrostatic interaction \cite{Sharp:1990a}.
\\

PB model is a one of the most used implicit solvent model in biological modeling \cite{Gilson:1988,Gilson:1987,David:2000}. It is a multiscale model that models the solution in two scales, the atomic modeling of the solute charges and the continuum modeling of the solvent ion distribution based on the Boltzmann distribution assumption. It is widely used in solvation free energy, binding free energy calculation in computational chemistry/biophysics community, the model can also be used for the molecular dynamics (MD) simulation, in which the force on the solute molecule due to the solvent effects is based on the PB calculation \cite{geng2011multiscale,Nina:1999,Gilson:1993,LuoRay:2013force}.
Despite the large amount of important potential applications of the PB model, it depends on the accurate solution to the PB model which yields the electrostatic field potential, many chemical and biological quantities can be derived from whom. The analytical solution to the PB model is rarely exists for the real solute molecules, the accurate and efficient numerical method for solving the PB model is important in making the model useful in the practical application. It directly promote the understanding of the structure, function, stability, and dynamics of solvated molecules.
\\

An accurate and efficient PB solver should address several issues. First issue in front of the development of the PB solver is the description of the solute molecular conformation structure. There are three main models used for modeling the conformational structure of the solvated molecule, namely the van der Waals (VDW) surface, solvent accessible surface (SAS), and the solvent excluded surface (SES). It is generally believed that the SES gives state-of-the-art accurate modeling of the solvated molecule, especially in terms of the electrostatics modeling in the implicit solvent scenario, which consistent with the explicit solvent simulation results \cite{LuoRayPB-VS-Explicitsolvent:2006, LuoRayPB-VS-Explicitsolvent:2010}. The second issue in developing numerical method for solving the PB model is the treatment of the singular charges in the source term of the governing equation, three typical treatments existing, the first one is projecting the charge to the grid points \cite{Jo:2008, Yu:2007}; the second one is treating the singular charge by the Green's function \cite{DuanChen:2011a, Geng:2007a, Chern:2003}; the third approach is named induced surface charge which project the induced charge to the molecular surface \cite{Chern:2003,LuoRayPBSA:2012,Delphi:2012}.
\\

Before developing the numerical method for the discretization the PBE, the solvated molecular structure should be defined. Many efforts have been paid to developing an accurate and robust SES software. Among those works, Connolly first invented a practical algorithm for implementing and visualizing of the SES \cite{SESConnolly:1983,Connolly:1983}. Later Sanner and his coworkers proposed a more robust way to calculate the reduced SES, which is named MSMS surface \cite{Sanner:1996}, the surface is able to treat the self intersecting surface in the SES. MSMS surface provides the surface area and volume calculation, also the triangulated surface is provided for further scientific computing purposes. Most recently, we have developed a Eulerian representation of the SES, the surface is named ESES \cite{ESES:2015}, the new surface is completely analytical which implements the SES and immersed to the Cartesian mesh through computing the intersections and normal directions of the surface and mesh lines. The ESES surface provides options for both finite element and finite difference calculations. Furthermore, the ESES surface can be utilized for the surface area, volume calculation, molecular cavities and loops detection. In parallel, in the past decades, many efforts also been devoted to approximating the SES, for instance, in the current widely used PB software package Delphi, an approximated SES was proposed, namely, smooth numerical surface (SNS) \cite{Rocchia:2002}, which provides a very good approximation of the SES for the electrostatics analysis of the biomolecules. Similar manner was adopted in another popular PB software in the Amber PBSA suites\cite{LuoRayPBSA:2012}.
\\

Due to the great importance of the numerical PB solver, in the past decades, many numerical method and PB solvers are developed by many different research groups. The numerical method for solving the PB model are generally classified into three categories, the finite difference method (FDM), the finite element method (FEM), and the boundary element method (BEM). The FDM is the most used discretization strategy among the PB solver community. Typical PB solvers based on the FDM formalism are Amber PBSA \cite{PBSA:2008,PBSA:2009}, APBS\cite{Holst:2000,Holst:2000}, Delphi\cite{Rocchia:2002,Delphi:2012}, CHARMM PBEQ\cite{Jo:2008}. The representative FEM category PB solver is the AFMPB \cite{LuBenzhuo:2013} which combines the Fast multipole method (FMM) with the adaptive FEM for solving the PB model. The BEM which use the potential theory to transform the PB model defined in the 3D domain to the molecular surface \cite{Juffer:1991}, and the PB model is solved on the triangulated molecular surface, a typical software applies this method is TABI-PB \cite{Geng2013:1,Geng2013:2,Geng2013:3}.
\\

In the MIB method, the interface conditions system have to be enlarged, so that the local complex geometry can be handled by eliminating the derivatives that are hard to be interpolated. According to the practical implementation, for the second order scheme, only two more interface conditions are needed. The important fact we used in developing the MIB method for solving the PB model is that, if a function is continuous across the interface, then so is its derivatives along the tangential directions.

In the past decade, we have dedicated ourselves to designing accurate and robust numerical schemes for solving   elliptic interface problems. In 2006, the matched interface and boundary (MIB) method \cite{Yu:2007c,Yu:2007a,Zhou:2006c,Zhou:2006d} was proposed, motivated by many practical needs, such as  optical molecular imaging \cite{DuanChen:2010b}, nano-electronic devices \cite{DuanChen:2010b},  vibration analysis of plates \cite{SNYu:2009}, wave propagation \cite{SZhao:2010a,SZhao:2008a},   aerodynamics   \cite{YCZhou:2012a}, elasticity \cite{BaoWang:2015c, BaoWang:2015d} and electrostatic potential in proteins \cite{Zhou:2008b,Yu:2007,Geng:2007a, DuanChen:2011a}. One important feature of the MIB method is its extension of computational domains by using the so called fictitious values, an strategy developed in our earlier methods for handling boundaries \cite{JCPWei:1999,Wei:2002f} and interfaces \cite{Zhao:2004}. As such, standard  central finite difference schemes can be employed to discretize differential operators as if there were no interface. Another distinct feature of the MIB method is to repeatedly enforce only the lowest order jump conditions to achieve higher order convergence, which is of critical importance for the robustness of the method to deal with arbitrarily complex interface geometries. Higher-order jump conditions must involve higher order derivatives and/or cross derivatives. Therefore, to approximate higher order derivatives and/or cross derivatives, one must utilize larger stencils, which is unstable for constructing high-order interface schemes and cumbersome for complex interface geometries.
Finally,  based on high order Lagrange polynomials, the MIB method is of arbitrarily high order in principle. For example, MIB schemes up to 16th order accurate have been constructed for simple interface geometries 1D and 2D domains \cite{Zhao:2004,Zhou:2006c}, and sixth-order accurate MIB schemes have been developed for complex interfaces in 2D \cite{Zhou:2006c} and 3D domains \cite{Yu:2007a,Zhou:2006c}. Recently we have constructed an adaptively  deformed  mesh  based MIB \cite{KLXia:2012a}  and a Galerkin formulation of MIB \cite{KLXia:2014e,KLXia:2014f} to improve MIB's capability of solving realistic problems. A comparison of the GFM, IIM and MIB methods can be found in in Refs. \cite{Zhou:2006c,Zhou:2006d}.
\\

Besides the accuracy of the PB solver, another important issue in the practical usage of the PB solver is its efficiency, the realistic biological system usually contains thousands to millions of atoms, the practical PB solver have to be fast enough to meet the requirements of simulating large biomolecules. There are mainly three types of approaches for speeding up the PB solver, the first category is the utilization of the out of core computing technique to carry out the computing in parallel and distributive manner, the second category is to accelerating the numerical discretization method in solving the PB equation. Another fundamental approach is to improve the accuracy and reducing the grid spacing influence of the PB solver itself, such that is accurate enough to implements the PB solver at very coarse grid.
\\

In our earlier work, we have introduced the MIB method and Green's function technique for treating the complex geometry and singular charges, respectively, in solving the PB model \cite{Geng:2007a,Yu:2007}. The method is generally tested to be of second order convergence in solving the PB equation, and provides better solvation free energy calculation than the classical PBEQ and APBS software for a testing set with 24 molecules. For the improvements, in this work, we firstly introduce our recently developed Eulerian represented analytical SES to represented the solvated molecules to avoid the error due to the approximation of the SES, it is tested that with the increasing of the surface density, the solvation free energy calculated on the MSMS surface will converge to ESES \cite{ESES:2015}. Second, we utilized a new strategy for the solvation free energy evaluation, in the previous work, the trilinear interpolation is adopted for the solvation free energy calculation, which is not appropriate due to the fact that the electrostatics potential in the solvent domain will be used which incurs the inaccuracy. Third, stabilization procedure (includes treating complex geometry and handle exception) is added to our earlier version package thus our software to solve the PB equation more robust and stable. The new software is tested to be highly accurate and robustness by both the Kirkwood solution and more than 1000 biomolecules. The current PB software can provide reliable solvation analysis at grid size as large as 1.1 to 1.2 \AA \ with the same level of accuracy as that from 0.2 to 0.3 \AA.
\\

This paper is organized in the following way. In section \ref{model} we briefly review the PB model which is formulated as an elliptic interface problem mathematically. In section \ref{numerical-method} we present the algorithms that used in solving the PB model. Section \ref{numericalresults-solvation} is devote to presenting the numerical results on electrostatics solvation free energy calculation, the numerical results validate the accuracy and robustness of our PB solver. The paper ends up with concluding remarks.

\section{Theory and method} \label{model}

\subsection{Theory}

\paragraph{The Poisson-Boltzmann equation}

In this section, we will give a brief review of the PB model. The PB model is a multiscale model which describes the solute molecule with atomistic detail while the solvent as a dielectric continuum \cite{Honig:1995a,Sharp:1990a}. Consider an open domain $\Omega\in \mathbb{R}^3=\Omega_m\bigcup\Gamma\bigcup\Omega_s$, where $\Omega_m$ is the solute domain that enclosed by the surface formed by the biomolecule, while $\Omega_s$ is the solvent domain, $\Gamma$ is the surface of the biomolecule that separates the solute and solvent domain. Let us  introduce a characteristic function $\chi ({\bf r}): \mathbb{R}^3 \rightarrow \mathbb{R} $ such that $\Omega_m  = \chi \Omega   $ and $\Omega_s  = (1-\chi) \Omega  $. Here, $\chi$ and $(1-\chi)$ are respectively the indicators  of the solute domain and the solvent domain.

The Poisson equation can be derived from the variation of the total free energy functional with respect to the electrostatic potential $\Phi$ \cite{Sharp:1990a, Gilson:1993, Geng:2011, JinPark:2015}
\begin{eqnarray}\label{eq24poisson}
-\nabla\cdot\left(\epsilon(\chi ) \nabla\Phi\right) =4\pi( \chi \rho_m + (1-\chi )\rho_s)
\end{eqnarray}
where $\epsilon(\chi )=(1-\chi )\epsilon_s+\chi \epsilon_m$ is the dielectric profile, which is $\epsilon_m$ in the solute domain and $\epsilon_s$ in the solvent domain. Here $\rho_m$ and $\rho_s$ are the charge densities of the solute and the solvent, respectively.

Due to the characteristic function $\chi$, the Poisson equation Eq. (\ref{eq24poisson}) is equivalent to
\begin{eqnarray}\label{eq24poisson2}
-\nabla\cdot\left(\epsilon_m \nabla\Phi\right) &=& 4\pi\rho_m,       \quad   {\bf r} \in \Omega_m \\ \label{eq24poisson2-2}
-\nabla\cdot\left(\epsilon_s  \nabla\Phi\right) &=& 4\pi\rho_s,         \quad {\bf r}   \in \Omega_s.
\end{eqnarray}
Mathematically, because electrostatic potential  $\Phi({\bf r})$  is defined on the whole computational domain ($\forall {\bf r} \in \Omega$), one has to impose
the following interface  conditions at the solvent-solute interface $\Gamma$ in order to make the Poisson equation Eq. (\ref{eq24poisson}) being well posed \cite{holst1994poisson,Geng:2007a,yu2007three}
\begin{eqnarray}\label{eq24poisson3}
[\Phi({\bf r})]    & = & 0, \quad  {\bf r}\in \Gamma; \\     \label{eq24poisson4}
[\epsilon({\bf r})\nabla \Phi({\bf r})] \cdot {\bf n} &=& 0, \quad {\bf r}\in \Gamma,
\end{eqnarray}
where $[ \cdot ]$ denotes the difference of the quantity ``$\cdot$" cross the interface $\Gamma$, ${\bf n}$ is  {interfacial norm direction} and
\begin{equation}
 \epsilon\left({\bf r} \right)=
\left\{\begin{array}{ll}\label{epsilonf}
  \epsilon_m, & {\bf r} \in\Omega_m; \\
  \epsilon_s, & {\bf r}\in\Omega_s.
\end{array}\right.
\end{equation}

For a second order PDE, its integration results in many arbitrary constants. Therefore, one must specify appropriate boundary condition to further make
the Poisson equation (\ref{eq24poisson}) being well posed. The variational derivation assumed the  far field boundary condition:
$$
\Phi(\infty)=0.
$$
However, in practical computation, the following Debye-H\"uckel boundary condition is employed:
\begin{equation} \label{DHBC}
\Phi({\bf r})=\sum_{i=1}^{N_m}\frac{Q_i}{\epsilon_s|{\bf r}-{\bf r}_i|}e^{-\bar{\kappa}|{\bf r}-{\bf r}_i|}, ~ {\bf r}\in \partial\Omega,
\end{equation}
where $\bar{\kappa}$ is the Debye-H\"uckel parameter \cite{ZhanChen:2010b}.

In Eq. (\ref{eq24poisson}), the solute charge density is given
 $$
    \rho_m({\bf r})=\sum_{i=1}^{N_m} Q_i\delta({\bf r}-{\bf r}_i), ~ {\bf r}\in  \Omega_m
$$
    where $Q_i$ is the partial charge of the $i$th atom, and $N_m$ is the number of charged atoms in the solute molecule.
Alternatively, one may compute the solute charge density directly by using the quantum mechanics \cite{ZhanChen:2011a}. In either treatment, the solute density gives an atomistic description of the solute. In contrast, the solvent charge density $\rho_s$ is described in a continuum manner and has the form
$$
\rho_s({\bf r})=\sum_{\alpha}^{N_c} \rho_{\alpha}({\bf r})q_{\alpha}, {\bf r}\in  \Omega_s
$$
where $ N_c$  represents the number of mobile ion species in the  solvent, and $\rho_{\alpha}({\bf r})$ and $q_{\alpha}$ are respectively the density and charge valence of the $\alpha$th  solvent species.  There are many ways to close Eq. (\ref{eq24poisson}) for $\rho_{\alpha}({\bf r})$. One way used in a non-equilibrium solvent,  is to let $\rho_{\alpha}({\bf r})$ be governed by another PDE, such as the Nernst-Planck equation   describing ion transport  \cite{JinPark:2015}. At equilibrium, one  assumes a Boltzmann distribution of the electrostatic energy for each ionic species in the Poisson--Boltzmann model,
 \begin{eqnarray}\label{eq24poissonB}
\rho_{\alpha}({\bf r}) =  \rho_{\alpha 0}e^{- q_{\alpha}\Phi({\bf r})\/k T},
\end{eqnarray}
where $\rho_{\alpha0}$ is the bulk density of $\alpha$th ion species, $k$ is the Boltzmann constant and $T$ is the temperature.

Equation (\ref{eq24poisson}), together with the Boltzmann distribution    of the electrostatic energy for solvent  ions  (\ref{eq24poissonB}),  interface conditions (\ref{eq24poisson3}) and boundary condition (\ref{DHBC}) is called the Poisson-Boltzmann equation (PBE).

\paragraph{Electrostatics solvation free energy}

The electrostatic solvation free energy, or reaction field energy, in the PB model can be represented as:
\begin{equation} \label{reaction-field-energy}
\Delta G_{\rm el}=\frac{1}{2}\sum_{i=1}^{N_m} Q_i \Phi_{\rm RF}({\bf r}_i),
\end{equation}
where  $\Phi_{\rm RF}({\bf r}_i)$ is the reaction field potential defined as $\Phi_{\rm RF}({\bf r}_i)\doteq \Phi({\bf r}_i)-\Phi_{\rm homo}({\bf r}_i)$ with $\Phi_{\rm homo}({\bf r}_i)$ being obtained by solving the PBE by switching  the dielectric constant $\epsilon_s$ in the solvent domain $\Omega_s$   to $\epsilon_m$.

In the rest of this paper, for the sake of simplicity, we only consider the case that there is no ion  in the solvent in describing the numerical method. The presence of ions in solvent does not affect our methods and accuracy of results.


\subsection{Numerical methods} \label{numerical-method}

This section presents a brief description of the MIB techniques for the solution of the PBE and the calculations of electrostatic solvation free energy.  Our MIB techniques consist of  three parts: 1) Interface setup and domain registration; 2)   Solution of the PBE, including the treatment of singular charges and the enforcement of the interface conditions; and 3) Calculation electrostatics solvation  free energy. These aspects are discussed in three subsections.

\subsubsection{Interface setup and domain registration}

In contrast to most PB solvers that do not strictly enforce interface conditions (\ref{eq24poisson3}), our MIBPB solver rigorously implement  interface conditions (\ref{eq24poisson3}) on the every intersecting points between the interface and mesh lines. As such, MIBPB techniques require to set up the interface explicitly and register the computational domain in terms of the solvent domain(s) and the solute domain.   We discuss the implementation of  the most popular and widely accepted surface,  solvent excluded surface (SES), in our computational techniques.  Our methodology can certainly be used for other smoother surfaces, such as Gaussian surfaces and minimal molecular surfaces \cite{Bates:2008}.

SES is defined as the trace of the probe that rolling around the atoms of the given solute molecule, which forms a cavity that cannot be penetrated by the solvent molecule. The SES is the most widely used surface definition \cite{Connolly:1983,Gilson:1988, Rocchia:2002}. With the SES description of the solvated molecular structure, the PB model can provide consistent results  in terms of electrostatic solvation free energies, electrostatic potentials of mean force of hydrogen-bonded and salt-bridged dimers calculation \cite{LuoRayPB-VS-Explicitsolvent:2006, LuoRayPB-VS-Explicitsolvent:2010, LuoRayPBSA:2012, Geng:2011}.

Explicit SESs are commonly expressed in two forms, namely, Lagrangian representation and Eulerian representation. The former is available from the  MSMS software package  \cite{Sanner:1996} and later can be obtained from the Eulerian  solvent excluded surface (ESES) software package \cite{ESES:2015}. To use MSMS surface for the finite difference based PB solvers, one needs to convert the triangular mesh surface into the surface description in the Cartesian domain. Whereas, ESES offers the Cartesian representation of SES directly.
Our MIBPB software we support both ESES and MSMS, and the numerical results in  Section \ref{numericalresults-solvation} indicates that MIBPB software can provides the same level of accuracy for the reaction field energy calculation on both surfaces.

\paragraph{MSMS surface in MIBPB}
The MSMS surface  is an efficient way for building the SES based on the reduced surface \cite{Sanner:1996}. It is the first SES computer program that can handle the geometric singularities  arising from  self-intersecting surfaces. The MSMS surface generation contains four steps for computing triangulated SES \cite{Sanner:1996}:
  a) Compute the reduced surface of a molecule;
 b) Build an analytical representation of the SES from the reduced surface;
 c) Remove all the self-intersecting parts from the analytical SES built above; and
 d)  Produce the triangulation of the reduced surface.
In our MIBPB software,  the MSMS surface in triangulated mesh has to be  embedded in the Cartesian mesh. To this end, we have developed projection algorithm for the Lagrangian to Eulerian transformation (L2ET)  \cite{YCZhouThesis}.
This transformation involves three main steps:
\begin{itemize}
\item Classify all the Cartesian grid points as either inside or outside the SES based on the ray-tracing, in which the discrete Jordan curve lemma is utilized \cite{Edelsbrunner:2010};

\item Calculate all the intersection points between  the triangulated SES and the Cartesian mesh lines. 

\item Calculate the associated interface normal directions at all the intersection points.
\end{itemize}
The above information is necessary and sufficient for  our MIB method to rigorously enforce interface conditions Eq. (\ref{eq24poisson3}) in solving the Poisson-Boltzmann equation.

\paragraph{ESES surface in MIBPB}
Recently, we have developed the ESES software package   \cite{ESES:2015} which directly provides all the information that is required in the rigious enforcement of the interface conditions Eq. (\ref{eq24poisson3}). The ESES is also implemented in the MIBPB software as a surface option. The benefits of the ESES is twofold. First, it guarantees the robustness of surface generation. In contrast, the MSMS (or the Lagrangian to Eulerian transformation) may fail to produce a surface at certain surface  density (or mesh size)  for certain proteins.   Additionally,   ESES generates  SES analytically in the Eulerian representation which avoids the error due to the triangulation of SES and the transformation of representations.  The SES generation in ESES software is based on the algorithm proposed by Connolly \cite{Connolly:1983}. To  accelerate the classification procedure of the Cartesian grid points, basic analytical  knowledge of the geometry and the ray-tracing technique are combined. To  find all the intersection points of the Cartesian mesh lines and the SES, ESES represents the surface patches of SES by  quadratic or quartic polynomials, and then solve appropriate  algebraic equations to locate intersection coordinates.  Advanced numerical solvers and computer graphics techniques are employed in developing the ESES software. The comparison of the reaction field potential based on both ESES and MSMS demonstrates that with the increasing of the MSMS surface density, the reaction field energy calculated by using MSMS will converge to that obtained by using ESES \cite{ESES:2015}.

\subsubsection{Solution of the Poisson-Boltzmann equation}

\paragraph{Green's function treatment of singular charges}
In most finite difference based  PB (FDPB) solvers, the singular charges located at the atomic centers are projected to the neighboring grid points of the Cartesian mesh. The problem with this projection is that when the grid size is larger than the  radius of  a solute atom near the interface, the atomic charge may be projected to the solvent domain. This unreasonable projection directly leads to the accuracy reduction of FDPB solvers at a coarse mesh.

Chern et al. proposed a mathematical procedure to analytically solve the singular charges by Green's functions  \cite{Chern:2003}. They have demonstrated their method by considering a 2D Poisson equation. The first Green's function treatment of realistic SESs of biomolecules was developed by us in conjugation with our MIB technique for the enforcement interface conditions \cite{Geng:2007a}. Let us consider the Poisson (\ref{eq24poisson}) without the salt charge ($\rho_s=0$).  In Green's function  approach,
%
%
the electrostatic potential  $\Phi$ is decomposed into the sum of the singular part ($\bar{\Phi}$ ) and regular parts  ($\tilde{\Phi}$):
\begin{equation}
\label{decomposition}
\Phi=\tilde{\Phi}+\bar{\Phi}.
\end{equation}
The singular part $\bar{\Phi}$ is defined as
\begin{equation}
\label{singular-part}
\bar{\Phi}({\bf r})=\left\{
                         \begin{array}{ll}
                           \Phi^*({\bf r})+\Phi^0({\bf r}), & \ {\bf r}\in \Omega_m, \\
                           0, & \ {\bf r}\in \Omega_s.
                         \end{array}
                       \right.
\end{equation}
where $\Phi^*$ is the Green's function
 \begin{equation} \label{GreenF}
\Phi^*({\bf r})=\sum_{i=1}^{N_m} \frac{Q_i}{\epsilon_m|{\bf r}-{\bf r}_i|}
 \end{equation}
that solves the Poisson equation with singular charges. Here
$\Phi^0$ is a harmonic function in the solute domain $\Omega_m$, which is obtained via solving the following boundary value problem (BVP):
\begin{equation} \label{BVP}
\left\{
  \begin{array}{ll}
    \nabla^2\Phi^0({\bf r})=0, & {\bf r}\in \Omega_m, \\
    \Phi^0({\bf r})=-\Phi^*({\bf r}), & {\bf r}\in\Gamma.
  \end{array}
\right.
\end{equation}

The regular part $\tilde{\Phi}$ satisfies the following homogenous  elliptic interface problem:
    \begin{equation}
    \label{eqnelliptical}
    -\nabla\cdot \left(\epsilon({\bf r})\nabla\tilde{\Phi}({\bf r}) \right)=0,
    \end{equation}
    with interface conditions:
    \begin{equation}
    \label{intface1}
    [\tilde{\Phi}]|_\Gamma=0, ~{\bf r}\in\Gamma
    \end{equation}
and
    \begin{equation}
    \label{intface2}
    [\epsilon \tilde{\Phi}_{\mathbf{n}}]|_\Gamma=-[\epsilon\bar{\Phi}_\mathbf{n}]|_\Gamma=\epsilon_m\nabla
(\Phi^*+\Phi^0)\cdot\mathbf{n}, ~{\bf r}\in\Gamma
    \end{equation}

    and the boundary condition:
    \begin{equation}
\label{RegBC}
\tilde{\Phi}({\bf r})=\sum_{i=1}^{N_m}\frac{Q_i}{\epsilon_s|{\bf r}-{\bf r}_i|},\   {\bf r}\in\partial\Omega.
\end{equation}

Obviously, the rigorous implementation of Green function approach requires the solution of the BVP Eq. (\ref{BVP}) with SES and the enforcement of interface conditions    Eq. (\ref{intface1}) and  Eq. (\ref{intface2}) in solving Eq. (\ref{eqnelliptical}).
These tasks are  technically as difficult as the solution of the original Poisson equation due to the  fact that   SES in the implicit solvent model setting  admits extremely complex geometry including geometric singularities, i.e., tips, cusps and self-intersecting surfaces.




\paragraph{Solution to the boundary value problem}

To obtain  the numerical solution $\Phi^0$ defined on the solute domain, we need to solve the BVP given by Eq. (\ref{BVP}). In our MIBPB software,  for the grid points at which the discretization does  not refer to any  grid point outside the solute domain, we direct employ the central FD (CFD) scheme; while for the other grid points, we need to interpolate the derivatives referred in the governing equation by both the grid values inside the SES and the boundary condition, where the nonuniform interpolation is carried out by using the Lagrangian polynomial based FD scheme \cite{Fornberg:1998}. For detail of the numerical method for solving the BVP,   the reader is referred to our earlier work \cite{Geng:2007a}.

\paragraph{Solution to the homogenous elliptic interface problem}
In this part, the MIB scheme for solving the elliptic interface problem Eqs. (\ref{eqnelliptical})-(\ref{RegBC}) is shortly reviewed. The computational domain $\Omega$ is discretized by the Cartesian mesh, with the grid size $dx, dy, dz$ along the $x, y, z$ direction, respectively. A grid point $(i, j, k) = (x_i, y_j, z_k)$  is said to be regular if all the 6 neighbor grid points and itself  referred  in a second-order CFD scheme are in the same subdomain, either $\Omega_s$ or $\Omega_m$. Otherwise $(i, j, k)$ is said to be irregular.
It is obvious that, at all the regular grid points the CFD scheme is applicable with the  second order convergence. However, at all the irregular grid points, due to the existence of the interface, the direct application of the CFD scheme  reduces the convergence order or accuracy. The main idea of the MIB scheme for solving the elliptic interface problem at the irregular grid point is to replace the function value at an irregular grid point in another subdomain by a fictitious value. The fictitious value is calculated via the rigorous enforcement of interface conditions \cite{Yu:2007,Geng:2007a, Zhou:2008b, DuanChen:2011a}.


\begin{figure}
\begin{center}
\begin{tabular}{cc}
\includegraphics[width=0.5\textwidth]{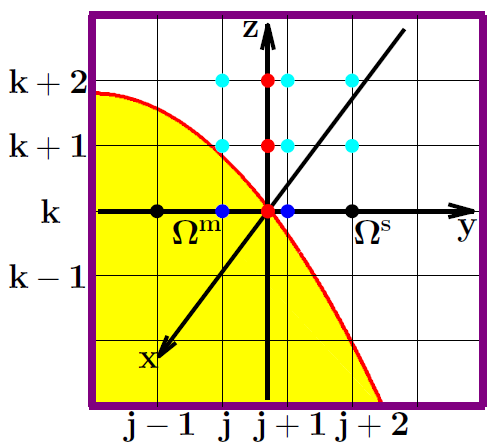}
\includegraphics[width=0.5\textwidth]{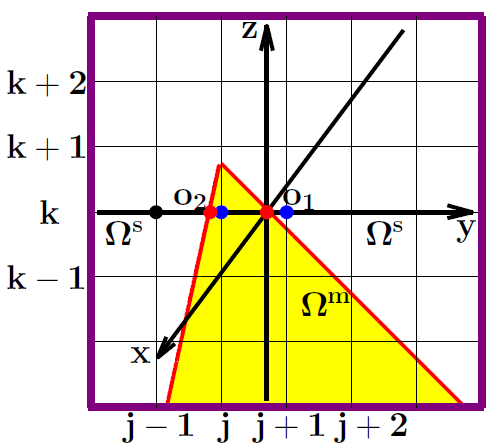}
\end{tabular}
\end{center}
\caption{Illustration of local interface geometries  around irregular grid point $(i, j, k)$. The left one is a locally smooth interface, while the right one is sharp interface. For the locally smooth interface, each time we solve two fictitious values at $(i, j, k)$ and $(i, j+1, k)$ along one direction simultaneously. For the sharp interface, each time we solve three fictitious values at $(i, j-1, k), (i, j, k)$ and $(i, j+1, k)$ along a direction involving two sets of interface conditions.
}
\label{local-geometry}
\end{figure}

According to the local geometry, the interface can be classified into two categories, locally smooth interface, and locally sharp interface. As depicted in Fig. \ref{local-geometry}, around a given irregular grid point $(i, j, k)$, the local interface is  smooth in the left chart and sharp in the right one. In both cases, we set up a local coordinate system at each intersection point between the interface and mesh lines. The relation between the general Cartesian mesh and the local coordinate is sorted out. Then  interface conditions   Eq. (\ref{intface1}) and  Eq. (\ref{intface2}) are expressed in the local coordinate. Additionally, it is useful to have more interface conditions so that more fictitious values can be determined. To this end, one might differentiates  interface conditions  Eq. (\ref{intface2}) to generate high-order derivatives, which typically leads to a  wider grid  stencils  and numerical instability. In our MIB method,  we avoid the use of high order interface conditions. Therefore, we create two new interface equations by differentiating Eq. (\ref{intface1}) along two  tangential directions. As such, we have a total of four interface equations that can be used to determine four factious values near each intersection point. Note that in these four interface conditions, there are six first order derivatives,  i.e., $\frac{\partial \tilde{\Phi}_s}{\partial x}, \frac{\partial \tilde{\Phi}_s}{\partial y}, \frac{\partial \tilde{\Phi}_s}{\partial z},  \frac{\partial \tilde{\Phi}_m}{\partial x}, \frac{\partial \tilde{\Phi}_m}{\partial y}, \frac{\partial \tilde{\Phi}_m}{\partial z}$, where subscripts ``$s$'' and ``$m$'' denote the limiting values at the interface $\Gamma$ from $\Omega_s$  and $\Omega_m$ domains, respectively, that are to be numerically computed at the intersection point. An MIB strategy is to further reduce the 3D interface problem into three 1D interface problems so that we only  need to determine two fictitious values along each mesh  line,
which gives us the liberty to eliminate two numerical derivatives near the intersection point that are most difficult to compute. This approach is used for locally smooth interface shown in the left chart of   Fig. \ref{local-geometry}. For locally non-smooth interface depicted in the right chart of   Fig. \ref{local-geometry}, a variety of sharp interface schemes were developed in our earlier work   to ensure the second order accuracy of our MIB method \cite{Yu:2007c,Yu:2007a,  KLXia:2011}. For example, one may consider two sets of interface conditions simultaneously and determines three fictitious values.

By the above MIB schemes,  all the fictitious values are determined  and the interface conditions are enforced. In the final FD discretization of the Poisson equation, when a derivative involves grid points at  a different domain, the fictitious values, instead of the original electrostatic potential  values, are used in the different domain.  We achieve the second order accuracy for arbitrarily complex interface geometry due to SESs \cite{Yu:2007c, Yu:2007, Geng:2007a}.


\subsubsection{Electrostatic solvation  free energy calculation}

One of major applications of the PB model is the calculation of the electrostatic solvation free energy, which offers the basis of the electrostatic binding energy as well.  The  accuracy, consistence, robustness and reliability of estimating biomolecular electrostatic solvation free energies provides an important  criteria for evaluating the performance of  PB solvers.

\paragraph{Reaction field potential representation of the solvation free energy}
In the conventional implicit solvent theory, the reaction field potential is defined as the difference between the electrostatic potential in solvent and in vacuum, that is,
\begin{equation}
\label{rec-potential1}
\phi_{\rm rec}(\mathbf{r})=\phi_{\rm dielec}(\mathbf{r})-\phi_{\rm vac}(\mathbf{r}),
\end{equation}
where $\phi_{\rm rec}, \phi_{\rm dielec}$ and $\phi_{\rm vac}$ are the reaction field potential, electrostatic potential in solvent and vacuum, respectively. The vacuum electrostatic potential is calculated through setting the solvent dielectric constant the same as that in the solute domain.

The electrostatics solvation free energy, or reaction field energy, is defined by:
\begin{equation}
\label{elec_eng}
\Delta G_{\rm RF}=\frac{1}{2}\sum_{i=1}^{N_m} Q(\mathbf{r}_i)\phi_{\rm rec}(\mathbf{r}_i),
\end{equation}
where $Q(\mathbf{r}_i)$ are the charge at the position $\mathbf{r}_i$.

\paragraph{Approximate the atomic center reaction field potential}
According to the equation for evaluating the electrostatics solvation free energy, Eq. (\ref{elec_eng}), the electrostatics potential at the atomic centers are
needed. Whereas, according to our numerical scheme, only the electrostatic potential on the grid of the Cartesian mesh grids is obtained, we need to approximate the electrostatic potential at the atomic centers by that at the grid points. In this work, the trilinear interpolation scheme will be utilized for interpolating the electrostatics potential at the atomic centers.

For $\forall \mathbf{r}_i\doteq (x_0, y_0, z_0)$, suppose the closest grid to $\mathbf{r}_i$ inside the solute domain is indexed by $(i, j, k)$. The following $27$ grids will be utilized for the further interpolation:
\begin{equation}
\label{inter-pts}
\{(i+m, j+n, k+p)|m=-1, 0, 1; n=-1, 0, 1; p=-1, 0, 1 \}.
\end{equation}

The trilinear interpolation scheme is given in the following steps:

\begin{itemize}
\item Interpolate the values at the following $9$ points:
$$
\{(x_0, j+n, k+p)|n=-1, 0, 1; p=-1, 0, 1\},
$$
by the above $27$ grids $\{(i+m, j+n, k+p)|m=-1, 0, 1; n=-1, 0, 1; p=-1, 0, 1 \}$, due to the utilization of the uniform Cartesian mesh, we have:
$$
\phi_{\rm rec}|_{(x_0, j+n, k+p)}=w_{x, -1}\phi_{\rm rec}|_{(i-1, j+n, k+p)}+
w_{x, 0}\phi_{\rm rec}|_{(i, j+n, k+p)}+w_{x, 1}\phi_{\rm rec}|_{(i+1, j+n, k+p)},
$$
for $n=-1, 0, 1; p=-1, 0, 1$, where $w_{x, -1}, w_{x, 0}$ and $w_{x, 1}$ are the Lagrangian interpolation coefficients, it is easy to be obtained via Fornberg's method \cite{Fornberg:1998}.

\item Interpolate the values at the following $3$ points:
$$
\{(x_0, y_0, k+p)| p=-1, 0, 1\},
$$
by the nine points $\{(x_0, j+n, k+p)|n=-1, 0, 1; p=-1, 0, 1\}$, similar to the above scheme, we have:
$$
\phi_{\rm rec}|_{(x_0, y_0, k+p)}=w_{y, -1}\phi_{\rm rec}|_{(x_0, j-1, k+p)}+w_{y, 0}\phi_{\rm rec}|_{(x_0, j, k+p)}+w_{y, 1}\phi_{\rm rec}|_{(x_0, j+1, k+p)},
$$
for $p=-1, 0, 1$, here $w_{y, -1}, w_{y, 0}$ and $w_{y, 1}$ are the interpolation coefficients.

\item Interpolate the value at the atom center $(x_0, y_0, z_0)$ by the $3$ points $\{(x_0, y_0, k+p)|p=-1, 0, 1\}$, which is given by:
$$
\phi_{\rm rec}|_{(x_0, y_0, z_0)}=w_{z, -1}\phi_{\rm rec}|_{(x_0, y_0, k-1)}+
w_{z, 0}\phi_{\rm rec}|_{(x_0, y_0, k)}+w_{z, 1}\phi_{\rm rec}|_{(x_0, y_0, k+1)},
$$
$w_{z, -1}, w_{z, 0}$ and $w_{z, 1}$ are the interpolation coefficients.
\end{itemize}

In sum, the approximation of the reaction field potential at the atomic center $\mathbf{r}_i$ is given by:
\begin{equation}
\label{27-interp}
\phi_{\rm rec}|_{(x_0, y_0, z_0)}=\sum_{m=-1}^1\sum_{n=-1}^1\sum_{p=-1}^1 w_{x, m}w_{y, n}w_{z, p}\phi_{\rm rec}|_{(i+m, j+n, k+p)}.
\end{equation}

\paragraph{The extension of the reaction field potential}
In the Eq.(\ref{27-interp}), for the atomic centers that close to the boundary when the coarse grid employed, it is possible that some of this grids may in the solvent domain, if we directly use the reaction field potential computed from the previous PB solver, the accuracy in evaluating the solvation free energy will be reduced. Therefore, we need to extend the reaction field potential in the solute domain to the outside grids that referred in interpolating the reaction field potential at the atomic centers.

For a given grid $(i_1, j_1, k_1)$ belongs to the above $27$ grids, the schemes that used for extension of the reaction field potential at $(i_1, j_1, k_1)$, $\phi_{\rm rec}(i_1, j_1, k_1)$, are listed in the following according to its priority, the scheme employed should have as top priority as possible.

\begin{itemize}

\item \textbf{Top priority:} Use the sum of fictitious value and the extended solution of the boundary value problem (specifically for the Green's function treatment of the singular charge case) at the grid point $(i_1, j_1, k_1)$ as the extended reaction field potential at $(i_1, j_1, k_1)$.

\item \textbf{Middle priority:} Choose three consecutive grid points next to $(i_1, j_1, k_1)$ along a given direction in the solute domain to extrapolate the reaction field potential at $(i_1, j_1, k_1)$.

\item \textbf{Low priority:} Select two inside grids neighbor to $(i_1, j_1, k_1)$ , say, $(i_2, j_2, k_2)$ and $(i_3, j_3, k_3)$ and approximate the reaction field potential at $(i_1, j_1, k_1)$ by:
$$
\phi_{\rm rec}(i_1, j_1, k_1)=\phi_{\rm rec}(i_2, j_2, k_2)+\phi_{\rm rec}(i_3, j_3, k_3)-\phi_{\rm rec}(i, j, k).
$$
\end{itemize}

\section{Validation and Results} \label{numericalresults-solvation}

In this section, we demonstrate the accuracy, consistence, robustness and reliability of our MIBPB method for the solution of the Poisson Boltzmann model. We utilize both analytical benchmark tests of Born model and Kirkwood model \cite{Kirkwood:1934} , and realistic proteins and protein complexes. All of these test examples  have been considered in earlier works \cite{LiAnbang:2014, Geng:2007a}.
In all tests, the grid dependence of the energies are examined. Errors are   estimated either with respect to the analytical solution if it is available or the result at the finest grid if there is no analytical result.
\begin{equation}
\label{measure}
{\rm Relative  \ error }\doteq \frac{| \Delta G_{h}-\Delta G_{\rm finest\  grid} |}{| \Delta G_{\rm finest\  grid}|},
\end{equation}
where $\Delta G_{h}$ and $\Delta G_{\rm finest\  grid}$ are the electrostatics solvation free energies calculated at the grid size $h$ and the finest grid used.
The dielectric constants adopted in the numerical tests are
$$
\epsilon({\bf r})=\left\{
                       \begin{array}{ll}
                         1, & \ {\bf r}\in\Omega_m \\
                         80, &\  {\bf r}\in\Omega_s.
                       \end{array}
                     \right.
$$

\subsection{Electrostatic solvation free energies}

In this subsection, we consider a large number of numerical tests on Electrostatic solvation free energies. The tests on the analytical Born sphere model and Kirkwood model reflect the accuracy of our methodology. A large amount of tests on biomolecules demonstrates both the grid size independent and robust properties of the MIBPB solver. Furthermore, the grid size independent property is verified on types of SESs, namely ESESs and MSMSs.

\subsubsection{Analytical tests}
In this part, the MIBPB software is  tested over examples with exact solutions.  The Born model for a dielectric sphere with a single point charge located at the atomic center and the Kirkwood model for dielectric sphere with multiple point charges are utilized in our study. For both the Born model and Kirkwood model \cite{Kirkwood:1934} the analytical solution is available. In fact, the Born model is a special case of the Kirkwood model.

\begin{table}[!ht]
\centering
\caption{Electrostatic solvation free  energies (kcal/mol) of a dielectric sphere with different radii and the central unit positive charge calculated by MIBPB.}
\footnotesize{
\begin{tabular}{lllllllllllll}
\hline
& & & & & &   {\rm Grid size} & ($h$) & & & & & \\
\cline{2-12}
{\rm Atomic ~ radius} &1.1  &1.0  &0.9 &0.8 &0.7  &0.6  &0.5 &0.4 &0.3 &0.2 &0.1 & Exact \\
\hline
1.1 (H)  &-143.63&	-116.22&	-148.63&	-145.9&	-146.03&	 -148.84&	 -149.00&	-148.98&	-149.01&	-149.04&	 -149.05&	 -149.05 \\
1.3 (H)  &-125.75&	-118.24&	-126.03&	-123.07&	-125.75&	 -125.99&	 -126.03&	-126.07&	-126.10&	-126.11&	 -126.12&	 -126.12\\
1.359(H) &-120.36&	-117.25&	-120.61&	-119.68&	-120.21&	 -120.54&	 -120.59&	-120.61&	-120.63&	-120.64&	 -120.65&	 -120.65\\
1.4 (O)  &-116.88&	-103.73&	-117.10&	-116.89&	-116.97&	 -117.03&	 -117.07&	-117.08&	-117.10&	-117.11&	 -117.11&	 -117.11  \\
1.5 (O)  &-109.17&	-103.43&	-109.18&	-109.17&	-109.22&	 -109.25&	 -109.28&	-109.29&	-109.29&	-109.30&	 -109.31&	 -109.13\\
1.55 (N) &-105.68&	-101.74&	-105.61&	-105.68&	-105.70&	 -105.72&	 -105.74&	-105.76&	-105.77&	105.78&	 -105.78&	 -105.78\\
1.7 (C)  &-94.44&	-95.80&	-96.34&	-96.40&	-96.39&	-96.40&	-96.42&	 -96.43&	 -96.44&	-96.44&	-96.45&	-96.45\\
1.8 (S)  &-90.96&	-90.96&	-91.02&	-91.01&	-91.04&	-91.05&	-91.07&	 -91.07&	 -91.08&	-91.09&	-91.09&	-91.09\\
1.85 (P) &-88.52&	-88.52&	-88.57&	-88.57&	-88.58&	-88.59&	-88.61&	 -88.62&	 -88.62&	-88.62&	-88.63&	-88.63 \\
2.0 (S)  &-81.86  &-81.95 &-81.90 &-81.92 &-81.95 &-81.96 &-81.98 &-81.97 &-81.98 &-81.98 &-81.98 &-81.98  \\
\hline
\end{tabular}
}
\label{Sphere_1_chg_mibpb}
\end{table}

\paragraph{Born model}

First, let us consider the dielectric sphere with a single unit charge at the atomic center center   and having different radii. Born theory states that, for a dielectric sphere with radius $R$ and center charge $Q$ placed in a solvent with dielectric constant $\epsilon_s$ and $\epsilon_m $ for the dielectric sphere, the electrostatic solvation free energy is:
\begin{equation}
\label{born}
\Delta G_{\rm el}=-\frac{Q^2}{2R}\left(\frac{1}{\epsilon_m} -\frac{1}{\epsilon_s}\right).
\end{equation}
In 2007, Born model was used to compare the performance of earlier versions of APBS, PBEQ and MIBPB  \cite {Geng:2007a}. Recently, Geng has used Born model to compare the accuracy of higher-order boundary integral equation approach with APBS and MIBPB \cite{Geng2013:2}. This comparison is indirect because the meshes of FD and boundary integral methods cannot be directed compared. Most recently,  Born model has been employed by Li to compare with the performance of DelPhi and Amber PB \cite{LiAnbang:2014}.

Table \ref{Sphere_1_chg_mibpb} lists the reaction field energies calculated from grid sizes 0.1 to 1.1\ \AA~ and the exact value. We consider radii used in the Amber force field for various atoms.  When the radius is 1.1\ \AA ~ and the mesh size is also 1.1\ \AA~, there geometry is hardly described by the mesh. One should not expect an accurate energy estimation. However,   MIBPB does an amazing job. In general, the radius is increased, the results on coarsest mesh size of 1.1\ \AA~ becomes more and more accurate. For example, when the radius is 1.8\ \AA and larger, relative errors on all grid sizes are less than 1\%. Additionally, for all radii, the relative errors on the mesh size of 0.9\ \AA ~ or smaller are always less than 1\%.
This result suggests that MIBPB can deliver reliable solvation estimation on a coarse mesh of 0.9\ \AA. Of course, the spherical geometry   is too simple for us to draw a conclusion, unless the MIB method can completely eliminate the geometric effect.

\paragraph{Kirkwood model}
To further test the accuracy of the PB solver, especially the Green's function treatment of the singular charges,
in this part we further test on the dielectric spheres, while now there are multiple charges distributed in the dielectric sphere.
The analytical solution to this case is due to Kirkwood \cite{Kirkwood:1934}.

We consider the following five different distributions of point charges but their radii are all set to be 2\AA.
\begin{itemize}
\item Case 1. Two positive unit charges symmetrically placed at $(1, 0, 0)$ and $(-1, 0, 0)$.
\item Case 2. Two positive unit charges symmetrically placed at $(1, 0, 0)$ and $(-1, 0, 0)$, and
       two negative unit charges symmetrically placed at $(0, 1, 0)$ and $(0, -1, 0)$.
\item Case 3. Two positive unit charges symmetrically placed at $(1.2, 0, 0)$ and $(-1.2, 0, 0)$,
      and two negative unit charges symmetrically placed at $(0, 1.2, 0)$ and $(0, -1.2, 0)$.
\item Case 4. Six Positive unit charges placed at $ (0.4, 0.0, 0.0), (0.0, 0.8, 0.0), (0.0, 0.0, 1.2), (0.0, 0.0, -0.4), (-0.8, 0.0, 0.0)$ and $(0.0, -1.2, 0.0)$.

\item Case 5. Six Positive unit charges placed at $(0.2, 0.2, 0.2), (0.5, 0.5, 0.5), (0.8, 0.8, 0.8), (-0.2, 0.2, -0.2), (0.5, -0.5, 0.5)$ and $    (-0.8, -0.8, -0.8)$\footnote{There were typos in the original coordinates for this test case in Ref. \cite{Geng:2007a}}.
\end{itemize}
The first three cases were employed by Li to compare Delphi and Amber PB performance \cite{LiAnbang:2014}. The last two cases were used in our earlier work in 2007 \cite{Geng:2007a}, which also contains many test examples similar to the first three cases. These benchmark tests are more difficult than the Born model. We highly recommend PB methodology developers to consider these test cases.

\begin{table}[!ht]
\centering
\caption{Electrostatic solvation free energies (kcal/mol) of Kirkwood dielectric sphere with multiple charges calculated by MIBPB}
\begin{tabular}{lllllllllll}
\hline
  & Case 1   && Case 2   && Case 3     & & Case 4  &&  Case 5  \\
\cline{2-3}
\cline{4-5}
\cline{6-7}
\cline{8-9}
\cline{10-11}
Grid size & $\Delta G_{\rm el}$  &RE &  $\Delta G_{\rm el}$ &RE &  $\Delta G_{\rm el}$  &RE  &   $\Delta G_{\rm el}$  &RE &  $\Delta G_{\rm el}$ &RE\\
\hline
$1.1$   &-351.12   &0.39\%  &-63.61  &1.27\%     &-135.41 &0.00\%  &-2974.59 &0.49\%  &-3079.70  &1.43\% \\
$1.0$   &-352.54   &0.80\%  &-65.66  &4.53\%     &-152.32 &12.45\% &-2952.83 &1.22\%  &-3138.85  &0.47\% \\
$0.9$   &-348.84   &0.25\%  &-60.70  &3.36\%     &-131.88 &2.59\%  &-2993.72 &0.15\%  &-3099.03  &1.80\% \\
$0.8$   &-351.20   &0.42\%  &-64.13  &2.10\%     &-137.33 &1.42\%  &-2996.27 &0.23\%  &-3110.38  &0.44\% \\
$0.7$   &-350.56   &0.24\%  &-63.24  &0.68\%     &-135.17 &0.16\%  &-2995.02 &0.19\%  &-3103.84  &0.65\% \\
$0.6$   &-350.68   &0.27\%  &-63.77  &1.52\%     &-137.34 &1.43\%  &-2994.73 &0.18\%  &-3118.62  &0.18\% \\
$0.5$   &-350.39   &0.19\%  &-63.50  &1.09\%     &-136.64 &0.92\%  &-2986.62 &0.09\%  &-3124.43  &0.00\% \\
$0.4$   &-350.02   &0.08\%  &-63.10  &0.46\%     &-136.27 &0.64\%  &-2991.70 &0.08\%  &-3119.22  &0.15\% \\
$0.3$   &-349.81   &0.02\%  &-61.95  &1.36\%     &-136.20 &0.59\%  &-2991.01 &0.06\%  &-3120.24  &0.13\% \\
$0.2$   &-349.72   &0.00\%  &-62.90  &0.14\%     &-135.79 &0.29\%  &-2989.89 &0.02\%  &-3122.89  &0.04\% \\
$0.1$   &-349.64   &0.02\%  &-62.81  &0.00\%     &-135.40 &0.00\%  &-2989.54 &0.00\%  &-3123.90  &0.01\% \\
\hline
Exact   &-349.73   &  &-62.81  &                 &-135.40 &        &-2989.30 &        &-3124.30  & \\
\hline
\end{tabular}
\label{MIBPBmulti_charge}
\end{table}

Table \ref{MIBPBmulti_charge} lists electrostatic solvation free energies of the MIBPB method for the above multiple charge tests on a set of mesh sizes range from 1.1 to 0.1\AA. For Case 1, all MIBPB relative errors are less than 1\%. For Case 2, MIBPB is smaller than 5\% on all grid sizes. Case 3 is relatively difficult because charges are located more close to the interface. For the last two cases, errors are bound by 1.5\% on all grid sizes.


\paragraph{Robust on different SES}
In this part, we  employ 25 biomolecules and data set  used in Li's  work  \cite{LiAnbang:2014} as the test set for comparing the performance of the Delphi \cite{Rocchia:2002} and Amber PBSA \cite{LuoRayPBSA:2012} PB software, the data was originally downloaded from the protein data bank. All the HETATM records in PDB files are removed by using MMTSB toolset, and the AMBER99SB force field was employed.  The structures were further optimized by the AMBER software \cite{LiAnbang:2014}.
The PDB IDs for this set of proteins are: 1ajj,  1ptq, 1vjw, 1bor, 1fxd, 1sh1, 1hpt, 1fca, 1bpi, 1r69, 1bbl, 1vii, 2erl, 451c, 2pde, 1cbn, 1frd, 1uxc, 1mbg, 1neq, 1a2s, 1svr, 1o7b, 1a63, and 1a7m.

\begin{table}[!ht]
\centering
\caption{Electrostatic solvation free energies (kcal/mol) calculated via the MIBPB software with MSMS surface at different grid sizes.}
\footnotesize{
\begin{tabular}{lllllllllll}
\hline
 & & & & &   {\rm Grid size} & ($h$) & & & & \\
\cline{2-11}
{\rm PDB ID}  &1.1  &1.0  &0.9 &0.8 &0.7  &0.6  &0.5 &0.4 &0.3 &0.2  \\
\hline
1ajj&-986.26 &-984.17 &-986.38 &-986.89 &-985.82& -986.68 &-986.64 &-986.85& -987.09 &-987.15  \\
1ptq &-723.95 & -722.99& -724.50 &-725.34 &-724.54  &-725.51 &-725.18 &-725.49& -725.64& -725.68  \\
1vjw &-1120.64 &-1117.29& -1119.03& -1119.37 &-1120.61 &-1121.60 &-1121.98 &-1121.52 &-1121.95 &-1122.12  \\
1bor&-773.51 &-774.38 &-777.01 &-776.91 &-778.85 &-778.04 &-778.33& -778.65 &-778.66& -778.7 \\
1fxd&-2424.97 &-2433.53 &-2432.13 &-2433.24 &-2434.06 &-2435.02 &-2434.64 &-2435.11 &-2435.14 &-2435.38  \\
1sh1 &-568.12 &-570.69 &-570.74 &-571.80 &-573.34&  -573.98 &-573.48 &-573.88 &-574.05 &-574.32   \\
1hpt &-662.19 &-666.16& -664.79 &-666.10 &-667.18 &-666.89& -667.57 &-668.17& -668.25 &-668.45  \\
1fca&-1175.78& -1181.58& -1181.77 &-1180.31& -1181.93& -1181.71& -1181.72 &-1182.04& -1182.10& -1182.11  \\
1bpi&-1150.37 &-1148.62 &-1149.68 &-1150.61 &-1151.08& -1151.53& -1152.48 &-1152.38 &-1152.63 &-1152.82  \\
1r69 &-945.34 &-938.93 &-940.63 &-940.05 &-942.28& -942.08 &-941.88& -941.98 &-942.15 &-942.23  \\
1bbl&-846.48 &-846.80 &-847.00& -846.90& -846.64& -847.33& -847.45& -847.08& -847.14 &-847.10 \\
1vii &-683.20 &-677.26 &-679.27 &-680.61 &-679.94 &-680.31& -681.19 &-681.28 &-681.40 &-681.55  \\
2erl &-919.91 &-915.41& -918.49 &-918.92 &-919.00& -919.68& -919.84 &-920.21& -920.33 &-920.57  \\
451c &-848.53 &-843.01 &-849.20 &-847.92 &-847.92& -847.63 &-847.85 &-848.21 &-848.03 &-848.07 \\
2pde &-798.26 &-799.19& -798.87 &-799.86 &-800.03  &-798.81& -799.39 &-799.44 &-799.22 &-799.26  \\
1cbn&-298.67& -301.97 &-303.85 &-303.68 &-304.11 &-304.46 &-304.41 &-304.79& -304.81 &-304.92   \\
1frd&-2561.35 &-2556.63& -2557.71& -2558.70& -2557.06& -2559.33 & -2560.09 & -2560.74 &-2561.27 &-2561.72 \\
1uxc &-917.39 &-914.88 &-918.19 &-919.73 &-920.18 &-921.01& -921.32&  -921.66 &-921.75& -922.02   \\
1mbg &-1286.82 &-1281.16 &-1283.11 &-1285.37 &-1285.74 &-1285.61& -1285.67& -1286.08& -1286.23 &-1286.45   \\
1neq &-1474.17 &-1471.80 &-1470.82 &-1472.88 &-1473.94 &-1474.13 &-1474.74& -1475.49 & -1475.62 &-1475.63 \\
1a2s &-1846.33& -1842.61& -1849.40& -1847.77& -1848.00 &-1847.38& -1848.25& -1848.48& -1848.86& -1849.09  \\
1svr &-1319.62 &-1320.81 &-1320.06 &-1323.10 &-1322.45 &-1323.89 &-1324.12 &-1324.76& -1324.68& -1324.99  \\
1o7b &-1738.56 &-1747.28 &-1743.36& -1743.60 &-1743.93 &-1745.48 &-1746.05 &-1746.18& -1746.77 & -1747.12  \\
1a63 &-1940.89& -1940.05 &-1936.96 &-1939.81 &-1941.67 &-1942.27 &-1942.61 &-1943.36 &-1944.14 &-1944.55  \\
1a7m&-2028.72& -2033.44& -2032.74& -2029.40& -2031.87& -2033.64& -2033.31& -2033.66& -2034.26& -2034.32   \\
\hline
\end{tabular}
}
\label{MIBPBmsms}
\end{table}

\begin{table}[!ht]
\centering
\caption{Electrostatic solvation free energies (kcal/mol) calculated via the MIBPB with ESES surface at different grid sizes.}
\footnotesize{
\begin{tabular}{lllllllllll}
\hline
 & & & & &   {\rm Grid size} & ($h$) & & & & \\
\cline{2-11}
{\rm PDB ID}  &1.1  &1.0  &0.9 &0.8 &0.7  &0.6  &0.5 &0.4 &0.3 &0.2  \\
\hline
1ajj &-996.11     & -995.63    & -993.17  &   -993.03  & -991.34  &   -988.16  &  -987.20 & -986.99 & -987.12 &   -987.03\\
1ptq &-720.68  &    -719.62    & -720.68  &   -720.76  & -721.48   &  -721.63  &  -721.71  &-722.04 & -721.95  &  -721.97\\
1vjw &-1122.67  &   -1123.35 &   -1124.63   & -1127.13 & -1124.10  &  -1125.55 &  -1125.15& -1125.20 &-1125.06 &  -1125.22\\
1bor &-771.50   &   -773.85    & -774.17 &    -772.50  & -774.31    & -773.74   & -774.04 & -774.37  &-774.44	 &  -774.45\\
1fxd & -2420.46   &  -2425.07 &   -2423.90 &   -2426.91 & -2427.53   & -2427.33  & -2428.19 &-2428.53& -2428.50  & -2429.05\\
1sh1 &-564.58   &   -564.16   &  -569.12   &  -569.23  & -569.74   &  -569.50  &  -570.31  &-570.59 & -571.01   & -571.20\\
1hpt &-661.72   &   -656.56  &   -661.03  &   -661.16  & -663.54  &   -663.66  &  -664.00  &-663.96 & -664.10   & -664.15\\
1fca &-1190.83    & -1193.85  &  -1191.26 &   -1190.80  &-1191.67  &  -1187.04  & -1188.83 &-1188.04 &-1187.68  & -1187.26\\
1bpi &-1139.46   &  -1145.12 &   -1147.15  &  -1145.35  &-1146.96  &  -1147.57 &  -1147.99 &-1148.14& -1148.28&   -1148.36\\
1r69 &-932.98  &    -937.44  &   -934.49  &   -936.89  & -937.27  &   -937.32  &  -937.10  &-937.44  &-937.59   & -937.66\\
1bbl &-838.68     & -839.89    & -841.26    & -842.17  & -843.02   &  -842.79  &  -842.00 & -842.80 & -842.73  &  -842.65\\
1vii &-672.12    &  -674.33 &    -674.94   &  -677.13 &  -676.71  &   -677.78 &   -677.69 & -678.05 & -678.42 &   -678.64\\
2erl &-914.31   &   -919.35   &  -918.70 &    -918.65 &  -917.98  &   -917.68 &   -918.68 & -918.60 & -918.79 &  -918.87\\
451c &-839.89   &   -838.39  &   -842.24  &   -841.94  & -842.56  &   -842.44  &  -842.56 & -843.26 & -843.78 &   -843.96\\
2pde &-809.06    &  -814.63  &   -812.91  &   -810.62 &  -814.13    & -813.32  &  -812.27 & -813.04& -813.59  &  -813.42\\
1cbn &-300.65     & -299.73    & -301.79  &   -302.88   &-301.94   &  -302.24  &  -302.37  &-302.50 & -302.60 &  -302.67\\
1frd &-2546.19   &  -2549.42 &   -2551.50  &  -2553.10  &-2554.68   & -2555.89  & -2556.43 &-2557.30 &-2557.69 &  -2558.10\\
1uxc &-909.77    &  -908.98  &   -912.19   &  -915.23  & -914.40   &  -916.34 &   -916.68 & -916.93 & -917.18 &   -917.44\\
1mbg &-1277.10   &  -1278.07  &  -1277.80 &   -1278.14 & -1280.52   & -1281.24  & -1281.42 &-1281.74 &-1281.86   &-1282.06\\
1neq &-1459.44    & -1464.56    &-1467.13   & -1467.57  &-1468.99  &  -1469.33 &  -1469.55 &-1469.64 &-1470.11  & -1470.12\\
1a2s &-1834.76     &-1838.39   & -1841.63   & -1842.01  &-1841.72  &  -1842.06  & -1842.67 &-1843.19& -1843.19  & -1843.61\\
1svr &-1309.47   &  -1310.67 &   -1314.35   & -1316.24 & -1317.25  &  -1317.10  & -1317.75& -1317.89& -1318.42  & -1318.37\\
1o7b &-1734.47   &  -1726.36  &  -1739.98  &  -1741.12&  -1739.71  &  -1740.59  & -1742.22 &-1742.57 &-1743.25 &  -1743.39\\
1a63 &-1937.29    & -1942.42    &-1940.78  &  -1943.18  &-1941.26  &  -1942.13  & -1941.64 &-1942.12 &-1942.27  & -1942.67\\
1a7m &-2038.64     &-2032.48    &-2034.48   & -2034.91  &-2038.38   & -2040.39  & -2037.91 &-2038.62& -2038.73  & -2038.26\\
\hline
\end{tabular}
}
\label{MIBPBeses}
\end{table}

Table \ref{MIBPBmsms} shows the electrostatic solvation free energies (kcal/mol)  calculated by using the MIBPB method  with SESs generated by the MSMS software at density 10. Similarly, Table \ref{MIBPBeses} shows the the electrostatic solvation free energies (kcal/mol)  calculated by the MIBPB software with SESs generated by the ESES software. First, there are some minor  discrepancies (less than 1\%) between the electrostatic solvation free energies computed by two SESs.  MSMS accuracy depends on its triangular mesh density and its results converge to those ESES at a very high triangular mesh density, say 100 triangles per \AA$^2$. Additionally, for a given surface, MIBPB is able to deliver very consistent results.
The results listed in the tables \ref{MIBPBmsms} and \ref{MIBPBeses} confirm the grid-size independence property of the MIBPB software.

\begin{figure}[!ht]
\small
\centering
\includegraphics[width=10cm,height=8cm]{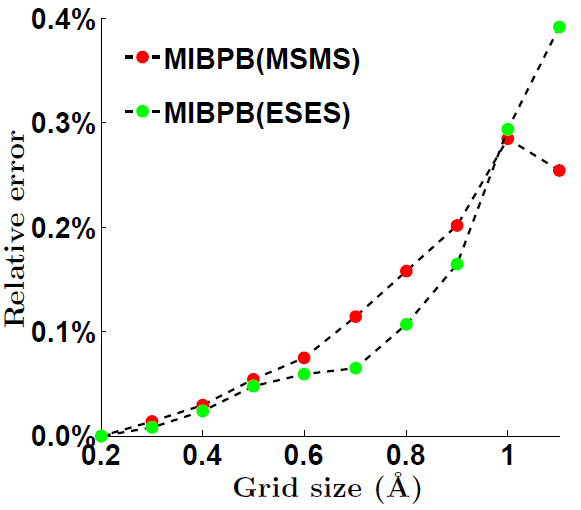}
\caption{The relative errors of the  electrostatic solvation free energies compared  to the results calculated at 0.2\AA~computed by the MIBPB method on surfaces generated by  ESES and MSMS averaged over   25 proteins.}
\label{rel-dev}
\end{figure}

Figure \ref{rel-dev} depicts the relative errors of the  electrostatic solvation free energies calculation by MIBPB software on surfaces generated by both ESES and MSMS. These results demonstrate that  MIBPB software on both surfaces are of  the same level the accuracy. Their relative errors are remarkably less than 0.4\% when the mesh size is refined from 1.1\AA ~ to 0.2\AA.  Therefore, if one's goal is 1\% relative error, one can just use a mesh size as coarse as 1.1\AA\  in MIBPB based solvation analysis.

\paragraph{Amber PB test set}

In this part, the MIBPB method will be tested on a much larger test set, namely, the Amber PBSA test set. This test set contains two parts, the nucleic acid and protein molecules, and has a total of 937 biomolecules. It provides the biomolecule structures, the corresponding charge and atomic radius parametrization. The data is available from \url{http://code.google.com/p/pbsa/downloads/detail?name=nucleicacidtest.tar.bz2} and \url{http://code.google.com/p/pbsa/downloads/detail?name=proteintest.tgz}.
 From our test, the MSMS software  cannot generate the SESs for the whole test set, whereas the ESES software can generate the surface successfully for all molecules. Therefore, we only report the MIBPB results of the electrostatic solvation free energies calculated  with ESESs.

\begin{figure}[!ht]
\small
\centering
\includegraphics[width=10cm,height=8cm]{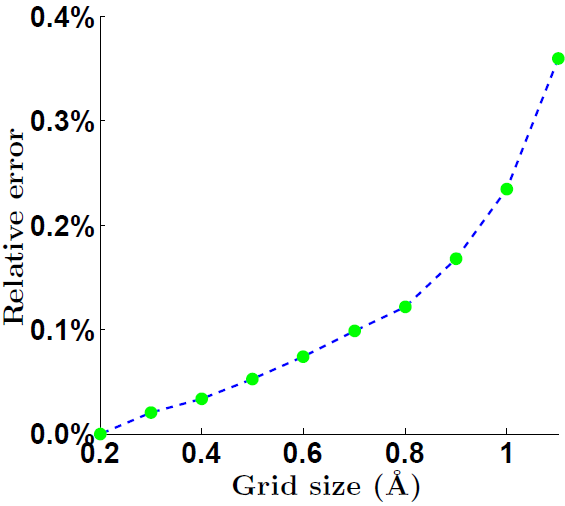}
\caption{The relative error of the electrostatics solvation free energies  with respect to those at grid size 0.2 \AA\ averaged over  937 biomolecules.}
\label{rel-dev-937}
\end{figure}

Figure \ref{rel-dev-937} shows the averaged relative errors of the electrostatics solvation free energies with respect to those at grid size 0.2 \AA \ for the 937 biomolecules. The errors  converge  monotonically with the refinement of the grid size, and the error level is consistent with that of the above 25 test set. Note that even at the grid size 1.1 \AA\  the averaged relative error is less than 0.4\%, which further indicates the grid size independence of our PB solver.

\label{conclusion}
\section{Concluding Remarks}
In this paper, we present a grid size almost independent PB solver, MIBPB software, it makes the accurate and coarse grid PB software possible. The main theme of the MIBPB solver is rigorously treating four details that referred in the PB model. First, the molecular surface definition, SES, is analytically implemented in the Cartesian grid for the purpose of the discretization of the PBE, which is different from the molecular surface used in most of the current existed PB solvers, most of which utilized the approximated SES instead of the exact SES. Second, the singular charges are treated by the Green's function instead of conventional method that project the singular charges to the closest eight grid points, the projection methods usually yields the charges be projected to the grid in the solvent domain when the coarse grid is employed, the Green's function treatment makes the coarse grid treating of the singular charges possible. Third, the interface conditions arise in the PBE are treated rigorously through the interface conditions matching, in the literature, there are some other methods that can treat these conditions with simple geometry while our method can handle geometry with arbitrary complex geometry. Fourth, the reaction field potential extension utilized for the evaluation of the reaction field energy, this posterior treatment avoids the accuracy reduction due to the usage of the reaction field potential in solvent.
\\

Due to the mathematical rigorously treatment and the utilization of second order convergence scheme, the MIBPB solver converge very fast. Which leads to the grid size almost independent property of the PB solver, this property was verified by a large amount of test cases includes both analytical tests and the real biomolecule tests, all the tests demonstrate that the current PB solver is grid size independent. Through the test of the stability with respect to different surfaces, the reduced molecular surface (MSMS surface) and the analytical solvent exclusive surface (ESES), the current PB solver's grid independent property is shown to be independent from the molecular surface used.
\\

In sum, this work developed a grid size independent PB solver, which makes the coarse grid PB solver becomes possible and indirectly speed up the PB solver.

\section*{Acknowledgments}
 This work was supported in part by NSF grant CCF-0936830, NIH grant R01GM-090208 and MSU Competitive Discretionary
Funding Program grant 91-4600. And we thanks to the ICER high performance computing center for the computing resources. Bao Wang thanks Dr. Geng Weihua, Dr. Li Anbang, Zhou Weijuan, and Ozsarfati Metin for lots of stimulating discussion and kindness help.

\end{document}